\documentclass{article}
\usepackage{amsmath,amsfonts,epsf,chicago,pb-diagram,lamsarrow,pb-lams,theorem}

\newcommand{\rf}[1]{(\expandafter\ref{#1})}
\newcommand{\ct}[1]{\citeANP{#1}~[\citeyearNP{#1}]}

\newcommand{\cttt}[3]{\citeANP{#1}~[\citeyearNP{#1}, \citeyearNP{#2}, \citeyearNP{#3}]}
\newcommand{\lb}[1]{\expandafter\label{#1}}

\newcommand\mfk\mathfrak
\newcommand\mcl\mathcal
\newcommand\mbb\mathbb
\newcommand\mtl\mathit
\newcommand\mbf\mathbf
\newcommand\onm\operatorname

\newcommand{\Hsp}{H^{\mbox{\scriptsize sp}}}
\newcommand{\omegasp}{\omega^{\mbox{\scriptsize sp}}}
\newcommand{\Jsp}{J^{\mbox{\scriptsize sp}}}
\newcommand{\Psp}{P^{\mbox{\scriptsize sp}}}

\newcommand{\Hpl}{H^{\!\mbox{\scriptsize pl}}}
\newcommand{\omegapl}{\omega^{\mbox{\scriptsize pl}}}
\newcommand{\Jpl}{J^{\mbox{\scriptsize pl}}}
\newcommand{\Ppl}{P^{\mbox{\scriptsize pl}}}

{\theorembodyfont{\rmfamily}

}

\begin{document}\allowdisplaybreaks
\title{Dynamics of Perturbed Relative Equilibria of  
Point Vortices on the Sphere or Plane}
\author{G. W. Patrick\\
Department of Mathematics and Statistics\\
University of Saskatchewan\\Saskatoon, Saskatchewan, 
S7N~5E6\\Canada
}
\date{June 1999}
\maketitle

\begin{abstract}
Stable assemblages  of localized vortices exist which have
particle-like properties, such as mass, and which can interact with
one another when they closely approach. In this article I calculate
the mass of these localized states and numerically investigate some
aspects of their interactions.
\end{abstract}
\section*{Introduction}

The system of $N$~point vortices on a sphere is a Hamiltonian system
with symmetry~$\mtl{SO}(3)$. Steadily rotating arrangements of
vortices are the relative equilibria of these systems, and the
nonlinear stability of these relative equilibria may be analyzed as a
matter of routine, at least in principle. If stable then such an
arrangement of vortices will, under sufficiently small perturbation,
approximately maintain its ``shape'', so that its subsequent motion
will be nearly its initial condition up to the action of the symmetry
group~$\mtl{SO}(3)$. However, if the momentum of the relative
equilibrium (i.e. its center of vorticity) is zero then the position
of such a relative equilibrium on the sphere is not stable under
perturbation to nonzero momentum, as has been established
in~\cttt{PatrickGW-1992.1}{PatrickGW-1995.1}{PatrickGW-1998.1}.

In this article I show that there are arbitrarily localized relative
equilibria of $N$~vortices with zero momentum, and that these relative
equilibria are formally stable for $N=3$ and $N=4$. Under small
perturbation to nonzero momentum, they move on the sphere,
to first order in the momentum perturbation, as a free particles move
under the influence of a magnetic monopole. The effect of the
monopole diminishes as the relative equilibria become more
localized. Since a given momentum perturbation implies a particular
velocity of motion of each relative equilibrium as a whole, they have
mass, and using a perturbation theory of~\ct{PatrickGW-1995.1}, I find
a formula for that mass.

For $N=4$, when two such relative equilibria are placed on the sphere,
they will move independently until they approach within a distance
comparable to their size, whereupon an interaction occurs. For
example, one might collide with its opposite (obtained by reversing
the first's component vortex strengths), the result being
``radiation'' in the form of dipole pairs, and generally opposites
repel and likes attract. Here in this article I numerically determine
some general features of these interactions, which are the subject of
ongoing investigations.

Over small regions of the sphere, the point vortex system on the
sphere becomes the point vortex system on the plane, and since the
above relative equilibria may be arbitrarily localized, one expects
such in the planar case as well. Towards the end of this article I
show this in indeed the case. The planar system of $N$ point vortices
will be a worthy motivating example for extensions of the stability
theory of Hamiltonian systems with symmetry to the case of
non-equivariant momentum mappings and/or noncompact groups.

\section{Context and Notations}
I briefly summarize here the basic elements of the system of $N$ point
vortices on the sphere and on the plane. For more details,
see~\ct{KidambiRNewtonPK-1998.1},~\ct{PekarskySMarsdenJE-1998.1}, and
the many references therein.

The system of $N$ point vortices of strengths $\Gamma_i\in\mbb R$ on
the sphere $S^2$ of radius $R$ is the Hamiltonian system with phase
space $\Psp\equiv (S^2)^N$ and Hamiltonian
\begin{equation}\lb{1}
\Hsp(x_1,\cdots,x_N)\equiv\frac1{4\pi R^2}\sum_{m<n}\Gamma_n\Gamma_m
 \onm{ln}{l_{mn}}^2,
\end{equation}
where $x_n$ represents the location of the $n^{\mbox{\scriptsize th}}$
vortex and $l_{mn}^2=2(R^2-x_m\cdot x_n)$ is the chord length between
$x_n$ and $x_m$. The symplectic structure $\omegasp$ on $\Psp$ is a
direct sum of the symplectic structures on $S^2$ weighted by vortex
strengths; specifically
\begin{equation}\lb{2}
\omegasp\equiv\bigoplus_{n=1}^N\frac{\Gamma_n}R\omega_{S^2},
\end{equation}
where
\begin{equation}\lb{3}
\omega_{S^2}(x)(v,w)\equiv-\frac1{R^2}x\cdot (v\times w).
\end{equation}
Most significant from a geometric mechanics point of view is that this
system admits the symmetry of diagonal multiplication of $\mtl{SO}(3)$
on $\Psp$ with momentum map
\begin{equation}\lb{4}
\Jsp\equiv-\frac1R\sum_{n=1}^N\Gamma_nx_n.
\end{equation}
These notations are as in~\ct{PekarskySMarsdenJE-1998.1}.

In the system of $N$ vortices in the plane the $n^{\mbox{\scriptsize
th}}$ vortex has location $z_n=x_n+iy_n\in \Ppl\equiv (\mbb C^2)^N$, and
the Hamiltonian and symplectic form are
\begin{equation}\lb{5}
\Hpl\equiv-\frac1{4\pi}\sum_{m<n}\Gamma_n\Gamma_m\onm{ln} |z_n-z_m|^2,\quad
\omegapl\equiv\bigoplus_{n=1}^N\Gamma_n\omega_0,
\end{equation}
where generally I will use the notation
$\omega_0(a,b)\equiv-\onm{Im}(a\bar{b})$ for complex numbers
$a,b\in\mbb C$. This system admits the symmetry group
$\mtl{SE}(2)=\{(e^{i\theta},a)\}$ of Euclidean symmetries acting
diagonally on each factor $\mbb C$ of $\Ppl$ by $(e^{i\theta},a)\cdot
z\equiv e^{i\theta}z+a$, and a momentum mapping is
\begin{equation}\lb{7}
\Jpl\equiv-\sum_{n=1}^N\Gamma_n\left[\begin{array}{c}\frac12|z_n|^2\\iz_n
\end{array}\right],\end{equation}
where $\mfk{se}(2)^*$ is identified with
$\mfk{se}(2)=\{(\dot\theta,\dot a)\}\cong\mbb R^3$ by the standard inner
product of $\mbb R^3$. These notations are as in~\ct{LewisDRatiuTS-1996.1}.

When attention is confined to a sufficiently small region of the
sphere, the spherical system becomes the planar one. For example, using
on each factor the map
\begin{equation}\lb{200}
\phi(z)=(z,\sqrt{R^2-|z|^2})
\end{equation} 
to pull back the
sphere to the plane one obtains
\begin{equation}\lb{101}
\phi^*\omegasp=-\frac1{zR}\omegapl\approx-\frac1{R^2}\omegapl,
\quad \phi^*\Hsp\approx-\frac1{R^2}\Hpl.
\end{equation}
Thus the pulled back equations of motion are approximately the same as
the equations of motion of the  planar system, since the factors of
$-1/R^2$ here multiply both the symplectic form and the Hamiltonian.

As noted in~\ct{AdamsMRatiuTS-1988.1}, the momentum mapping of the
planar system is equivariant if and only if the total vortex strength
$\sum_n\Gamma_n$ vanishes. Indeed, the adjoint and coadjoint
actions of $\mtl{SE}(2)$ are
\begin{gather}\lb{8}
\onm{Ad}_{(e^{i\theta},a)}(\dot\theta,\dot a)=(\dot\theta,e^{i\theta}a
 -i\dot\theta a),\\
\onm{CoAd}_{(e^{i\theta},a)}(\mu,\nu)=\bigl(\mu-\omega_0(
e^{i\theta}\nu, a),e^{i\theta}\nu\bigr),\nonumber
\end{gather}
and the deviation of the momentum map $\Jpl$ from equivariance is the
cocycle $\sigma:\mtl{SE}(2)\rightarrow \mfk{se}(2)^*$ defined by
$\sigma\equiv \Jpl\bigl((e^{i\theta},a)\cdot
p\bigr)-\onm{Coad}_{(e^{i\theta},a)}\Jpl(p)$ (the evaluation point
$p\in P$ is irrelevant). For the case at hand the cocycle is
\begin{equation*}
\sigma(e^{i\theta},a)=-{\textstyle\bigl(\sum_n\Gamma_n\bigr)}
 \left[\begin{array}{c}\frac12|a|^2\\ia\end{array}\right],
\end{equation*}
and the derivative of this cocycle at the identity is the skew-symmetric
two form $\Sigma:\mfk{se}(2)^2\rightarrow\mbb R$ given by
\begin{equation}\lb{104}
\Sigma\bigl((\dot \theta_1,\dot a_1),(\dot \theta_2,\dot a_2)\bigr)
={\textstyle\bigl(\sum_n\Gamma_n\bigr)}\omega_0(\dot a_1,\dot a_2).
\end{equation}
Generalities on non-equivariant momentum mappings may be found
in~\ct{AbrahamRMarsdenJE-1978.1}; a prime fact is the momentum
commutation identity
\begin{equation}\lb{20}
\{J_\xi,J_\eta\}=J_{[\xi,\eta]}-\Sigma(\xi,\eta).
\end{equation}

Both the planar and spherical systems have simple closed form
solutions in the case $N=2$: for the sphere any two vortices evolve as
the action of the one parameter group with generator $\Jsp/2\pi
{l_{12}}^2$ while for the plane they evolve as the action of the one
parameter group with generator
$\dot\theta=(\Gamma_1+\Gamma_2)/2\pi|z_1-z_2|^2, \dot
a=-\Jpl/2\pi|z_1-z_2|^2$. Numerically the action of these one
parameter subgroups is easily computed, and since the Hamiltonians
$\Hsp$ and $\Hpl$ are sums of pairwise interactions, the full system
of $N$ vortices may be numerically integrated in a symplectic,
symmetry preserving and momentum preserving way using splitting
methods, as in~\ct{ChannellPJNeriFP-1993.1}.

\section{The Relative Equilibria}
To find the relative equilibria for the system of $N$ vortices in the
sphere one seeks $\xi_e\in\mfk{so}(3)\cong\mbb R^3$ and $p_e\in\Psp$
such that $d\Hsp(p_e)-dJ_{\xi_e}(p_e)=0$. Equivalently, using the
obvious extensions of $\Hsp$ and $\Jsp$ to $(\mbb R^3)^N$, one can
solve the equations
\begin{equation*}
\nabla_m\Hsp-\nabla_m\Jsp_{\xi_e}=\tilde\lambda \nabla
\left(\frac12|x_m|^2\right)=\tilde\lambda_mx_m,\quad m=1,\cdots,N,
\end{equation*}
where the $\tilde\lambda_m\in\mbb R$ are Lagrange
multipliers. Inserting~\rf{1} and~\rf{4} yields the $N$ equations
\begin{equation}\lb{10}
\sum_{n\ne m}\frac{\Gamma_n}{l_{nm}^2}=\lambda_mx_m+2\pi R\xi_e,\quad
\lambda_m\equiv-\frac{2\pi R^2}{\Gamma_m}\tilde\lambda_m,\quad m=1,\cdots,N.
\end{equation} 
General analytic solutions are not to be expected to these nontrivial
nonlinear equations.

I seek a manifold of solutions to~\rf{10} contained in the zero level
set of $\Jsp$, which in some limit is confined to arbitrarily small
regions of phase space, and which is formally stable near that limit.
I try a regular polygonal configuration where the first $N-1$ vortices
of equal strengths surround the (possibly different strength)
$N^{\mbox{\scriptsize th}}$ central vortex, $N\ge3$. I will denote the
strength of the central vortex by $\Gamma$. For convenience I locate
the central $N^{\mbox{\scriptsize th}}$~vortex at $R\mbf k$ where
$\mbf k\equiv(0,0,1)$, and the first vortex at
$(R\sin\alpha,0,R\cos\alpha)$, so that $\alpha$ is the opening
angle---the angle at the center of the sphere between any outer vortex and
the central vortex. The momentum is zero if and only if
\begin{equation*}
\Gamma_1(x_1+x_2+\cdots+x_{N-1})+\Gamma x_N=(N-1)\Gamma_1R\cos\alpha
\mbf k+\Gamma R\mbf k=0
\end{equation*}
which is equivalent to
\begin{equation}\lb{11}
\Gamma_1=-\frac{\Gamma }{(N-1)\cos\alpha},
\end{equation}
while equations~\rf{10} reduce to $\lambda_1=\cdots=\lambda_{N-1}$ and the
two equations
\begin{gather}
\lb{12}\Gamma_1X_0+\frac{R\Gamma }{l_{1N}^2}\mbf k=\lambda_1x_1+2
 \pi R\xi_e,\\
\lb{13}\frac{(N-1)\Gamma_1R\cos\alpha}{l_{1N}^2}\mbf k=\lambda_NR\mbf k
 +2\pi R\xi_e,
\end{gather}
where $\xi_e$ depends on $\Gamma$, $N$, $\alpha$ and $R$, and where
\begin{equation}\lb{14}
X_0\equiv\frac{x_2}{l_{21}^2}+\frac{x_3}{l_{31}^2}+\cdots+
\frac{x_{N-1}}{l_{N-1,1}^2}.
\end{equation}
Setting
\begin{gather*}
\mu_N\equiv\sum_{k=1}^{N-1}\frac1{e^{-2\pi ik/N}-1},\quad
\nu_N\equiv\sum_{k=1}^{N-1}\frac1{|e^{2\pi ik/N}-1|^2},
\end{gather*}
so that
\begin{equation*}
\mu_N+\nu_N=\sum_{k=1}^{N-1}\frac{e^{2\pi ik/N}}{|e^{2\pi ik/N}-1|^2},
\end{equation*}
one sees by symmetry that $\mu_N$ and $\nu_N$ are real, and also that
\begin{equation*}\mu_N=-\frac12(N-1),\quad \nu_N=\frac1{12}(N^2-1).
\end{equation*}
Then scaling $\mu_{N-1}$ and $\nu_{N-1}$ gives
\begin{equation}\lb{15}
X_0\cdot\mbf i=\frac{\mu_{N-1}+\nu_{N-1}}{R\sin\alpha},\quad 
X_0\cdot\mbf j=0,\quad
X_0\cdot\mbf k=\frac{\nu_{N-1}\cos\alpha}{R\sin^2\alpha}.
\end{equation}
Equation~\rf{13} implies $\xi_e$ is along $\mbf k$ and in~\rf{12}
it is clear that the $\mbf j$ component of both sides is zero, so
writing~\rf{12} and~\rf{13} in components gives three linear equations
in the three unknowns $\lambda_1$, $\lambda_N$, and
$\xi_e\cdot\mbf k$. These equations can be routinely solved to
obtain
\begin{gather}
\lambda_1=-\frac{\Gamma }{R^2\cos\alpha\sin^2\alpha}
 \frac{(N-2)(N-6)}{12(N-1)},\lb{16}\\
\lambda_N=-\frac{\Gamma }{R^2\sin^2\alpha}\left(\cos\alpha+\frac{N}{2(N-1)}
 \right),\lb{17}\\
\xi_e\cdot\mbf k=\frac{\Gamma }{4\pi R^2\sin^2\alpha}\left(\frac1{N-1}
 +\cos\alpha\right).\lb{18}
\end{gather}
In the limit $\alpha\rightarrow 0$, from~\rf{11}, the
total vortex strength
\begin{equation}\lb{*3}
{\textstyle\sum_n\Gamma_N}=\Gamma +(N-1)\Gamma_1=\Gamma \left(1-\frac1{\cos\alpha}\right)
\end{equation}
vanishes while the regular polygon of outer vortices collapses upon
the central vortex. Although these relative equilibria depend on
$\alpha$, $\Gamma$ and $N$, I will denote them below simply by $p_e$.

Formal stability means definiteness of the Hessian of the reduced
Hamiltonian at the relative equilibrium, and, by the energy-momentum,
method as in~\ct{MarsdenJE-1994.1}, this is equivalent to definiteness
of the function $\Hsp-J_{\xi_e}$ on a subspace tangent to the momentum
level set and complimentary to the tangent space to the group orbit.
For $N=3$ the reduced spaces are points and formal stability is
immediate. I have determined the stability $p_e$ for $N=4,5,6$ with the aid
of a symbolic manipulator. For $N=4$ the relative equilibria are
formally stable without conditions, and in particular formally stable
arbitrarily near the limit $\alpha\rightarrow0$. For $N=5$, but
numerically now, one has formal stability if and only if
$\alpha>1.951\approx111.8^o$ and for $N=6$ if and only if
$\alpha>2.245\approx128.6^o$. Presumably this pattern continues and
the relative equilibria are formally stable for $N>6$ if and only if
$\alpha$ is sufficiently near $\pi$.

As already mentioned, since the momenta of such relative equilibria
are zero, by~\ct{PatrickGW-1992.1}, formal stability implies stability
only modulo~$\mtl{SO}(3)$; dynamically, these relative equilibria,
when perturbed to nonzero  momentum, approximately maintain their
shape, which oscillates on a fast time scale, while they move around
the sphere on a slow time scale. The left of Figure~\rf{25} is the
result of a simulation and shows a typical motion of one such relative
equilibrium. For short, I will call a point  of phase space that
results  from perturbing one of  the  relative equilibria $p_e$ a
``preq''.

The motion of preq may be approximated as a direct
application of the perturbation theory of~\ct{PatrickGW-1995.1}. In
summary, the symplectic $\mtl{SO}(3)$ symmetry implies that the
linearization of the relative equilibrium has double eigenvalues and
hence can be expected to have a nonzero nilpotent part, say
$N_\alpha$. It can be shown that the image of $N_\alpha$ is contained
in the tangent space to the group orbit and the tangent space to the
momentum level set is contained in the kernel of $N_\alpha$.
Consequently there is a unique (it can be shown to be symmetric)
bilinear form (here with the same name) $N_\alpha$ on $\mfk{so}(3)^*$
such that the following diagram commutes
\begin{equation*}\begin{diagram}
\node{T_{p_e}\Psp}\arrow{e,t}{N_\alpha}\arrow{s,l}{dJ(p_e)}
 \node{T_{p_e}\Psp}\\
\node{\mfk{so}(3)^*}\arrow{e,t}{N_\alpha}
 \node{\mfk{so}(3)}\arrow{n,r}{\xi\mapsto\xi p_e}
\end{diagram}\end{equation*}
I have calculated the $N_\alpha$ in the case $N=3$ and $N=4$; they are
the diagonal matrices with entries $J_{N1},J_{N2},J_{N3}$ where
\begin{gather}
J_{31}=J_{32}\equiv\frac1{8\pi R^2\cos\alpha},\lb{31}\\
J_{33}\equiv\frac1{8\pi R^2\sin^4\alpha}(2\cos^3\alpha+3\cos^2\alpha
+2\cos\alpha-1),\lb{32}
\end{gather}
and
\begin{gather}
J_{41}=J_{42}\equiv
-\frac1{2\pi R^2}\frac{3\cos^2\alpha+\cos\alpha+2}{
(1-\cos\alpha)(9\cos^2\alpha+4\cos\alpha+3)},\lb{33}\\
J_{43}\equiv\frac1{12\pi R^2\sin^4\alpha}(3\cos^3\alpha+4\cos^2\alpha
+3\cos\alpha-2).\lb{34}
\end{gather}
The perturbation theory asserts that to first order the motion of the
preq  is that  of the \emph{drift  system}: the reduction, by the
\emph{normal form symmetry} of the right-hand action of the torus
generated by $\mbf k\in\mfk{so}(3)$, of the left invariant \emph{drift
Hamiltonian} $\pi\mapsto\pi\cdot\xi_e+\frac12\pi^t  N_\alpha\pi$  on
$T^*\mtl{SO}(3)$ with canonical symplectic form. That the drift Hamiltonian
is invariant under this toral action is predicted by the theory and is
evident by the equality $J_{N1}=J_{N2}$. This flavor of the theory is
subject to the caveat the the frequencies of the linearization of the
equilibrium are not in 1-1 resonance with the rotation frequencies of
the relative equilibrium itself; were there to be such a 1-1 resonance
the preq would acquire the  structure of a  magnetic moment, as
in~\ct{PatrickGW-1998.1}. For the case $N=3$ this is not an issue
since the reduced spaces are points and so there are no reduced
frequencies. For $N=4$ one verifies by calculating the linearization
that if the reduced frequencies are $\pm\omega_{\mbox{\scriptsize
red}}$ then
\begin{equation*}
\frac{|\xi_e|^2}{\omega_{\mbox{\scriptsize red}}^2}=1-
 3\frac{(9\cos^2\alpha+4\cos\alpha+3)\sin^2\alpha}
 {3\cos^2\alpha+2\cos\alpha+3}.
\end{equation*}
Since the second term does not vanish for $0<\alpha<\pi$, the group
frequencies $\pm\xi_e$ are never equal to the reduced frequencies and
the simpler perturbation theory suffices.

If the location of the preq on the sphere is denoted
by $y$, then the equations of motion for it (i.e. the
equations of motion of the drift system) turn out to be
\begin{equation}\lb{85}
m_\alpha\frac{d^2y}{dt^2}=-m_\alpha\frac{|v|^2}{R}\frac{y}{R}
+\sigma B\times v,
\quad v\equiv\frac{dy}{dt},\quad B\equiv\frac1{R^3}y,
\end{equation}
where $\sigma$ the momentum associated to the normal form symmetry and
\begin{equation*}
m_\alpha\equiv-\frac{1}{J_{N1}R^2}.
\end{equation*}
Thus, the drift system is the same as that of particle of mass
$m_\alpha$ and charge $\sigma$ moving on the sphere under the
influence of the magnetic monopole $B$. Given that one is perturbing
a relative equilibrium where the central vortex has been located at
$R\mbf k$, the initial location of the preq may be taken to be
$y(0)=R\mbf k$, while its initial velocity~$v(0)$ and the
charge~$\sigma$ may be obtained from the momentum
perturbation~$\Delta\mu$ by the equations
\begin{equation*}
v(0)=\frac{1}{m_\alpha R^2}y(0)\times\Delta\mu,\quad
\sigma=-\Delta\mu\cdot\mbf k.
\end{equation*}
By solving Equations~\rf{85}, one finds that the prediction of the
theory is that the preq will rotate about the perturbed momentum, say
$\Delta\mu$, with angular velocity $J_{N1}\Delta\mu$.

Specified momentum perturbations of the relative equilibria $p_e$ may
be accomplished by moving the central vortex from its original
location at $R\mbf k$ to $(R\sin\delta,0,R\cos\delta)$ while changing
the angle $\alpha$ to $\alpha+\Delta\alpha$, as follows. The momentum
of the perturbed configuration is easily verified to be
\begin{equation*}
\Jsp=-(n-1)\Gamma_1
\left[\begin{array}{c}0\\0\\\cos(\alpha+\Delta\alpha)\end{array}\right]
-\Gamma \left[\begin{array}{c}\sin\delta\\0\\\cos\delta-1
\end{array}\right],
\end{equation*}
so this is $\Delta\mu\equiv[\Delta\mu_1,0,\Delta\mu_3]$ if
\begin{gather}
\sin\delta=-\frac{\Delta\mu_1}{\Gamma },\lb{26}\\
\cos(\alpha+\Delta\alpha)-\cos\alpha=-\frac1{(N-1)\Gamma_1}\bigl(\Delta\mu_3
+\Gamma (\cos\delta-1)\bigr),\lb{27}
\end{gather}
and arbitrary directions in $\Delta\mu$ may be obtained by rotating
this. As these formulas are to be used to perturb a relative
equilibrium it is understood that $\Delta\mu_1$ and $\Delta\mu_3$ are
small; just how small depends on the validity of the approximation
that is the drift system. The momentum of the relative equilibrium,
being zero, does not itself provide a scale. Certainly, though, the
perturbation should not significantly affect the geometry of the
relative equilibrium, meaning it should not displace the vortices by
amounts comparable to the diameter $2\alpha$. In the using of~\rf{26}
and~\rf{27} one should therefore ensure
\begin{equation}\lb{28}
\left|\frac\delta\alpha\right|\ll1,\quad
\left|\frac{\Delta\alpha}{\alpha}\right|\ll1
\end{equation}
Using~\rf{26} and~\rf{27} I have simulated the $4$-vortex preq
corresponding to $\alpha=\pi/6$ on the unit sphere for momenta
$\Delta\mu$ small multiples of the vector $[2,0,3]$, and calculated
the angular velocities of the preq. The results are on the right of
Figure~\rf{25}. The drift rates fit well to a cubic with slope
$|J_{41}|$ at zero.{\renewcommand{\baselinestretch}{.8}\footnotesize
\begin{figure}
\setlength\unitlength{1in}
\centerline{\begin{picture}(4.5,2.1)
\put(0,.20){\epsfbox{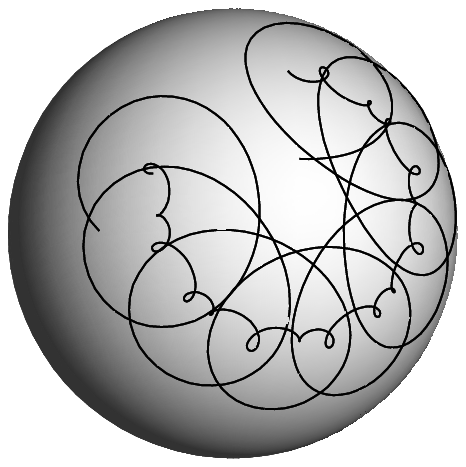}}
\put(2.2,.6){\epsfysize=1.5in\epsfbox{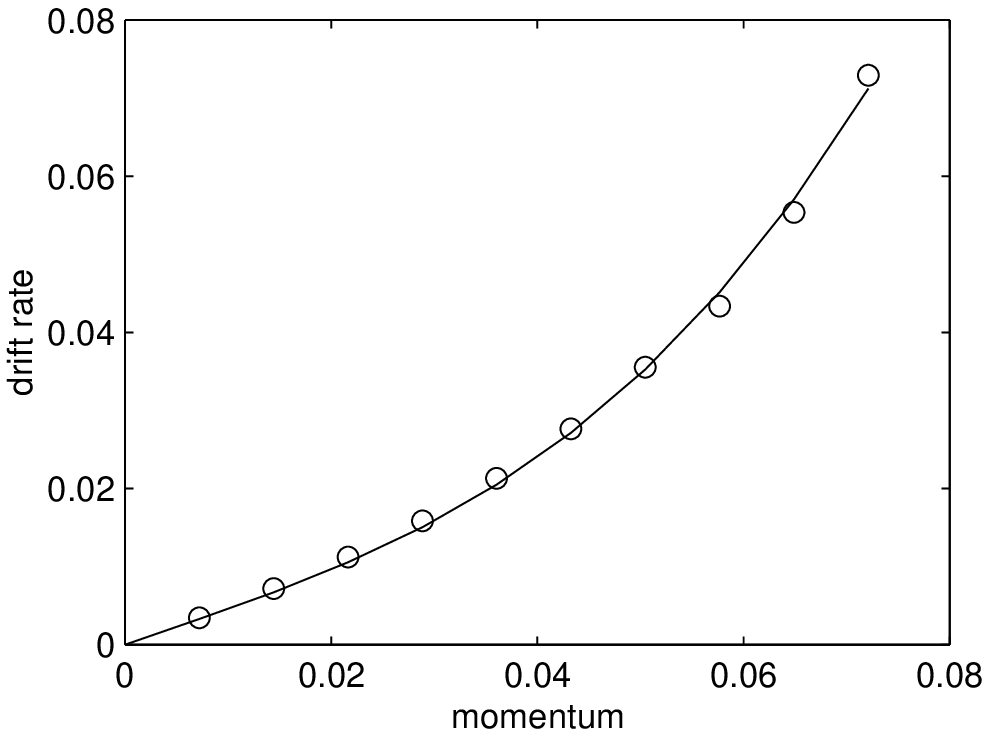}}
\put(2.2,0){\epsfxsize=2.25in\epsfbox{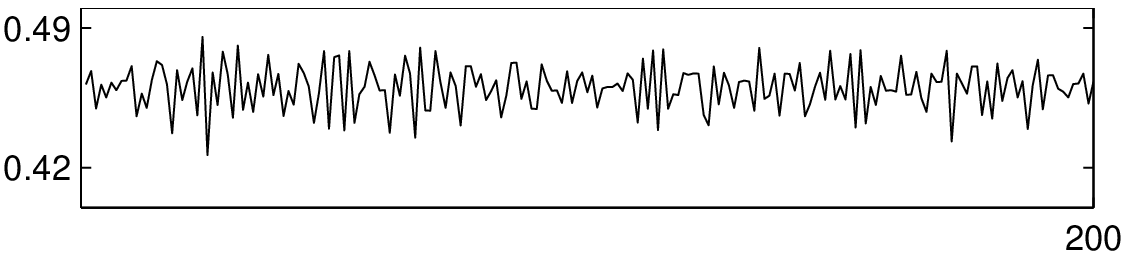}}
\end{picture}}
\caption{\lb{25} { \it Left: A $4$-vortex preq with opening angle
$\alpha=\pi/6$ subjected to a large perturbation subsequently moves on
the sphere in a circular path about the perturbed momentum, which
points towards you. The sense of rotation as is clockwise, the
reverse of the right hand sense obtained from the perturbed angular
momentum. Shown are the paths of the central vortex and one outer
vortex; the paths of the other two outer vortices are suppressed to
avoid cluttering the picture. Upper right: the drift rate vs
$|\Delta\mu|$; the curve shown is a cubic with slope at $0$ equal
to $J_{41}$.} Lower right: the angular velocity of a preq divided by
$|\Delta\mu|$, calculated over $200$ consecutive intervals of its path.
\normalsize}
\end{figure}}

An essential aspect of the preq is that their location is ill defined
as a concept. Suppose for example one assigns the location $R\mbf k$
to a preq corresponding to the state $p_0\in\Psp$ at some particular
time, and sometime later the state of the system is $p_1\in\Psp$. If
there is a group element $A\in\mtl{SO}(3)$ such that $Ap_0=p_1$ then
the location corresponding to $p_1$ is $RA\mbf k$, unequivocally.
However, for $N\ge 4$ the the reduced phase space for the vortex
system has dimension $2N-4\ge4$, so the flow on this reduced phase
space is usually at least as complicated as a toral flow with two
incommensurate frequencies. Thus the reduced flow may never repeat
itself, and there may never be an $A$ such that $Ap_0=p_1$ exactly.
For the simulations above the location of the preq was
decreed to be the average of the locations of its constituent
vortices; another possibility for example is the location of its
central vortex. The preq alters its shape on a fast time
scale compared with its overall motion, and this gives a statistical
character to the meaning of its ``location'' on the sphere. The
location of the preq becomes ever more exact and its character
ever less statistical as the perturbation falls to zero.

This problem is illustrated in the bottom right of Figure~\rf{25},
where path of a preq has been regularly sampled and
plotted are the values of $J_{41}$ from angular velocities obtained
from consecutive changes in the angle that the preq makes with
its initial condition. Were the preq to exactly follow the drift
system the plot would be constant, but instead one gets a variation
about the average value of $0.4602$. This average is just $.07\%$ off
the value of $J_{41}=.4599$ predicted by Equation~\rf{33}.

The issue of whether it is possible to define the location of a
preq may be put more deeply as follows. Let
the Marsden-Weinstein reduced space for the vortex system at its zero
momentum level be $\Psp_0$ with reduced Hamiltonian $\bar\Hsp$; the
relative equilibrium $p_e$ corresponds to the equilibrium, say $\bar
p_e$ of the reduced system. In~\ct{PatrickGW-1995.1} it is shown that
there is a symplectomorphism defined near the group orbit
$\mtl{SO}(3)\cdot p_e$ and onto a neighborhood of $\bar p_e$ times the
zero section of $T^*\mtl{SO}(3)$ such that the
Hamiltonian $\Hsp$ becomes
\begin{equation*}
\tilde\Hsp(x,\pi)\equiv\bar\Hsp(x)+\pi\cdot\xi_e+\frac12 N_\alpha\pi^2
 +\text{h.o.t.}
\end{equation*}
Thus the reduced degrees of freedom are linked to the overall motion
by the higher order terms depending on both $x$ and $\pi$, and if
these terms vanish then the location of the preq may be
unequivocally assigned by following the fiber of $T^*\mtl{SO}(3)$ to
its zero section. Conversely, the location of the preq is ill defined
in as much as one is obstructed in removing these higher order terms.

\section{Small Opening Angles}

With opening angles as large as $\alpha=\pi/6$, such as in
Figure~\rf{25}, the relative equilibrium $p_e$ does not have the appearance
of a localized particle. But one may choose $\alpha$ arbitrarily
small, and now I discuss some aspects of small $\alpha$ and report
some numerics for $\alpha=\pi/2048\approx.088\onm{deg}$.

Firstly, in the $\alpha\rightarrow0$ limit there is an essential
difference between the 3-vortex and 4-vortex preq: for
the 3-vortex preq, from~\rf{31},
\begin{equation*}
m_\alpha=-8\pi+O(\alpha^2),
\end{equation*}
while for the 4-vortex preq, from~\rf{33},
\begin{equation*}
m_\alpha=\frac83\pi\alpha^2+O(\alpha^4).
\end{equation*}
Since the mass of the 3-vortex preq does not fall with $\alpha$ while
that of the 4-vortex preq does, the 3-vortex preq is very heavy in
comparison to the 4-vortex preq for small $\alpha$. Thus, large
momentum perturbations are required to move a 3-vortex preq. As shall
be shown immediately, momentum perturbations must also fall with
$\alpha$, or else destroy the relative equilibrium by too grossly
perturbing it. Thus, for dynamical purposes the 3-vortex preq is
infinitely heavy in the $\alpha=0$ limit. This is confirmed upon
passage to the planar system, as will be seen in Section~(4).

The restrictions~\rf{28} in the $\alpha=0$ limit have a different
character as they affect $\Delta\mu_1$ vs. $\Delta\mu_3$. The effect
on $\Delta\mu_1$ is clear: directly from~\rf{26} and the first
of~\rf{28}
\begin{equation*}
|\Delta\mu_1|=|\Gamma \sin\delta|\ll O(\alpha).
\end{equation*}
However, by elementary manipulations,~\rf{27} becomes
\begin{equation*}
|\Delta\alpha|=|\cos^{-1}(z\cos\alpha)-\alpha|,
\end{equation*}
where, temporarily,
\begin{equation}\lb{29}
z\equiv\cos\delta+\frac{\Delta\mu_3}{\Gamma }.
\end{equation}
By considering the function $z\mapsto|\cos^{-1}(z\cos\alpha)-\alpha|$
for small $\alpha$ one verifies that
$|\cos^{-1}(z\cos\alpha)-\alpha|\ll\alpha$ is equivalent to
\begin{equation*}
0<\frac{z-1}{\frac1{\cos\alpha}-1}<<1\quad\mbox{or}\quad
0<\frac{1-z}{\frac{\cos2\alpha}{\cos\alpha}-1}<<1.
\end{equation*}
In any case this amounts to $|z-1|<O(\alpha^2)$ and since
$\delta\ll\alpha$ one gets by~\rf{29} that $\Delta\mu_3\ll
O(\alpha^2)$. The point is that for momentum perturbations
$\Delta\mu_1$, $\Delta\mu_3$ to be considered small the first must
fall as $\alpha$ while the second must fall as $\alpha^2$. Since
$|\Delta\mu_1|\gg|\Delta\mu_3|$ causes preq motion that is rotation about an
axis perpendicular to the direction of the preq, the
effect is that the preqs can have momentum $O(\alpha)$
only for motions $O(\alpha)$-close to great circle paths. To obtain
preq motion along smaller circular paths on the sphere one
must take $\Delta\mu_1\approx\Delta\mu_3=O(\alpha^2)$, which implies a
much smaller $O(\alpha^2)$ momentum.

Since the maximum reasonable momentum perturbations are $O(\alpha)$
and the 4-vortex mass is $O(\alpha^2)$, the maximum 4-vortex angular
velocity on the sphere is $O(1/\alpha)$. Thus greater velocities are
available to smaller opening angles; the velocity available to a
4-vortex preq compared with its diameter is $O(1/\alpha^2)$.

{\renewcommand{\baselinestretch}{.8}\footnotesize
\begin{figure}
\setlength\unitlength{1in}
\centerline{\begin{picture}(3.75,1.5)
\put(0,0){\epsfbox{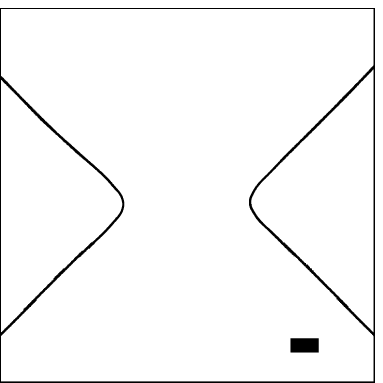}}
\put(2.25,0){\epsfbox{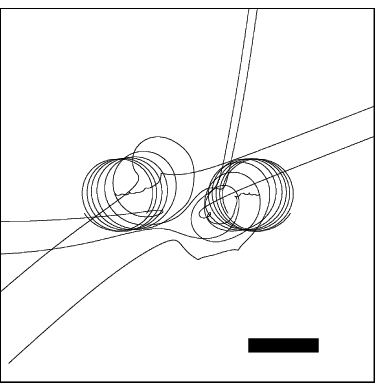}}
\end{picture}}
\caption{\lb{80}{\it Left: A preq incident from the lower
left approaches its anti-preq incident from the lower right; a
repulsive force is evident. Right: the same preq and
anti-preq collide head on and break up into dipole pairs. The
bars at the lower left indicate the diameter of the preqs.
\normalsize}}\end{figure}}

A  most interesting aspect of preqs  with small opening angles is that
more than one of them  may be positioned  on the sphere, all initially
far  apart  (compared with  their radii)  from one  another. Since the
vortex-vortex interaction  falls  as  the  the  vortex-vortex distance
increases, as  long as the separate preqs   remain separated they will
not strongly  interact with one  another, and they will  maintain their
separate  identities. However, as they  separately move  on the sphere
they may closely approach one another. When  two or more preqs closely
approach the above  perturbation theory becomes  inapplicable and they
undergo  an interaction, and may move   apart from one another largely
unchanged,  or may be partially or   completely destroyed. There is no
theory of these interactions  at  this  time  but the  main  features,
obtained  by numerical simulation  and illustrated  in Figures~\rf{80}
and~\rf{81}, are as follows:
\begin{itemize}
\item Two preqs interact when their distance falls to lengths
comparable with their diameters.
\item
Two preqs with vortices having the same sign are attracted to one
another while if the vortices have opposite signs they repel one
another. Thus identical preqs attract while a preq and its anti-preq
(obtained from the first by changing the signs of each vortex) repel.
\item 
The collision of two identical preqs, partially since those preqs
attract, tends to result in two groups of four vortices representing a
state of the vortex system far away from the relative equilibria
$p_e$. In this sense the 4-vortex preqs are usually destroyed in
same-sense vortex preq collisions. In such an interaction the two
4-vortex preqs usually exchange one of their outer vortices.
\item
A very energetic collision of a 4-vortex preq with its anti-preq
usually results in the destruction of both preqs into vortex dipole
pairs.
\item
The detailed results of an energetic or non-elastic 4-vortex preq
collisions are extremely sensitive to initial condition, while the
gross aspects of the collision (e.g. whether or not the collision is
elastic) are relatively more robust.
\end{itemize}
Many a happy hour may be spent watching the antics of 4-vortex preqs
colliding of the sphere and the nature of the 4-vortex interaction is
the subject of ongoing investigations.

{\renewcommand{\baselinestretch}{.8}\footnotesize
\begin{figure}
\setlength\unitlength{1in}
\centerline{\begin{picture}(3.75,1.5)
\put(0,0){\epsfbox{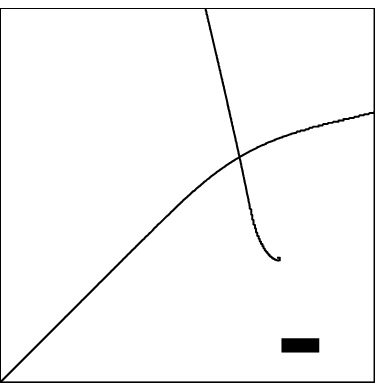}}
\put(2.25,0){\epsfbox{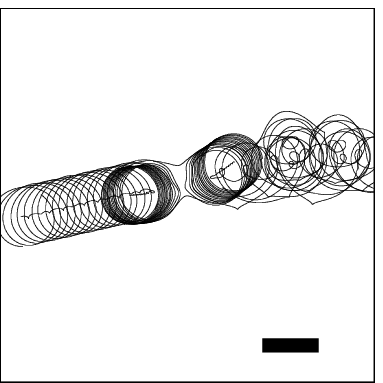}}
\end{picture}}
\caption{\lb{81}{\it Left: A preq incident from the lower left
approaches an identical preq initially at rest in the center; and
attractive force is evident. Right: a preq moving rightward and
slightly upward collides with an identical preq moving leftward and
slightly downward. The denser tracks at the center are the incident
preqs. The preq incident from the left remains and rebounds back the
way it came. The preq incident from the right has been nearly
destroyed and is a four vortex motion far from the original relative
equilibrium. Clearly visible in the interaction is the exchange of an
outer vortex. \normalsize}}\end{figure}}

\section{Transcription to the Planar System}

As already noted (see Equations~\rf{101}), when restricted to
sufficiently small regions of the sphere, the system of point vortices
on the sphere reverts to the system of point vortices on the
plane. This is pertinent for the study of small-opening-angle preq
interactions, since the interactions occur only as the preq closely
approach. Moreover, the planar case, being an example with noncompact
symmetry group, is currently  interesting, as the general focus of the 
Hamiltonian stability
literature has been on the compact case, while the noncompact
situation is known to have distinctive features, as shown
by~\ct{LeonardNEMarsdenJE-1997.1}. Consequently, I close this article
with a short transcription of the above the the planar case; the
results should be viewed as the $\alpha=0$ limit of the spherical
system.

The transcription is easy: through the pull back~\rf{200} of the sphere
to the plane, the relative equilibria $p_e$ with opening angle
$\alpha$ have radii $\tilde\alpha=R\sin\alpha$, and by~\rf{18} they
rotate at angular velocity
\begin{equation}
\dot\theta_\alpha=\frac\Gamma{4\pi R^2\sin^2\alpha}\left(\frac1{N-1}
+\cos\alpha\right)=\frac N{N-1}\frac\Gamma{4\pi\tilde\alpha^2}
+O(\tilde\alpha^4),
\end{equation}
while, from~\rf{11},
\begin{equation*}
\Gamma_1=\frac{-\Gamma}{(N-1)\cos\alpha}=\frac{-\Gamma}{N-1}
+O(\tilde\alpha^2).
\end{equation*}
Deleting the higher order terms and renaming $\tilde\alpha$ to
$\alpha$, one expects planar relative equilibria with radii $\alpha$,
central vortex strength $\Gamma$, outer vortex strength
$\Gamma/(N-1)$, and $\mfk{se}(2)$ generator $(\dot\theta_\alpha,0)$,
where
\begin{equation*}
\dot\theta_\alpha=\frac{N}{N-1}\frac\Gamma{4\pi\alpha^2},
\end{equation*}
and it is easily verified that these are indeed relative equilibria of
the planar system. Moreover, these relative equilibria occur for the
vortex strength parameters $\Gamma_n$ such that $\sum\Gamma_n=0$,
parameters for which, by~\rf{104}, the momentum map $\Jpl$ is
equivariant.

For my purpose, equivariance of the momentum map is important since, in
the general, equivariance guarantees, for a relative equilibrium $p_e$
with momentum $\mu_e$ and generator $\xi_e$, the
\emph{momentum-commutation relation} $\onm{coad}_{\xi_e}\mu_e=0$, and
this is crucial for the perturbation theory
of~\ct{PatrickGW-1995.1}. On the other hand, in absence of
equivariance one has the more complicated commutation relation
\begin{equation*}
\onm{coad}_{\xi_e}\mu_e=-\Sigma(\xi_e,\cdot)=-i_{\xi_e}\Sigma,
\end{equation*}
and the perturbation theory would have to be extended in a fundamental
way to cover the nonequivariant case. That not being necessary, one
calculates the nilpotent parts of the linearizations of the planar
system at the relative equilibria above, to obtain, for the 3-vortex
relative equilibria
\begin{equation*}
N_\alpha=\frac{1}{2\pi\alpha^4}\left(\begin{array}{ccc}
-\frac34&0&0\\0&0&0\\0&0&0\end{array}\right),
\end{equation*}
while for the 4-vortex relative equilibria the nilpotent part turns out
to be
\begin{equation*}
N_\alpha=\frac{1}{2\pi\alpha^2}\left(\begin{array}{ccc}
-\frac43&0&0\\0&\frac34&0\\0&0&\frac34\end{array}\right),
\end{equation*}
Going to the Lagrangian viewpoint, and comparing with the Lagrangian
of a 2-dimensional rigid body moving in the plane, one sees that the
3-vortex preq have infinite mass while the 4-vortex preq have inertia
coefficient $-\frac3{2\pi}\alpha^2$ and mass $\frac{8\pi}3\alpha^2$.

\footnotesize\frenchspacing

\end{document}